\documentclass[12pt,twoside]{article}
\usepackage[english]{babel}
\usepackage[latin1]{inputenc}
\usepackage{amsmath}
\usepackage{amssymb,amsfonts}
\usepackage{graphicx}
\usepackage{times,amssymb,amscd}

\newcommand{\bA}{\mathbf{A}}

\newcommand{\bE}{\mathbf{E}}

\newcommand{\bH}{\mathbf{H}}

\newcommand{\bL}{\mathbf{L}}

\newcommand{\bR}{\mathbf{R}}

\newcommand{\bS}{\mathbf{S}}
\newcommand{\bV}{\mathbf{V}}

\newcommand{\bs}{\mathbf{s}}
\newcommand{\ba}{\mathbf{a}}
\newcommand{\bb}{\mathbf{b}}

\newcommand{\be}{\mathbf{e}}

\newcommand{\br}{\mathbf{r}}
\newcommand{\bx}{\mathbf{x}}
\newcommand{\by}{\mathbf{y}}

\newcommand{\BV}{\boldsymbol{V}}

\newcommand{\Be}{\boldsymbol{e}}
\newcommand{\Bu}{\boldsymbol{u}}
\newcommand{\Bv}{\boldsymbol{v}}

\newcommand{\cF}{\mathcal{F}}

\newcommand{\cP}{\mathcal{P}}
\newcommand{\cS}{\mathcal{S}}
\newcommand{\cT}{\mathcal{T}}

\newcommand{\cB}{\mathcal{B}}

\newcommand{\cM}{\mathcal{M}}

\newcommand{\EUC}{\bE^3}

\newcommand{\SXR}{\bS^2\!\times\!\bR}
\newcommand{\HXR}{\bH^2\!\times\!\bR}
\newcommand{\SLR}{\widetilde{\bS\bL_2\bR}}
\newcommand{\NIL}{\mathbf{Nil}}
\newcommand{\SOL}{\mathbf{Sol}}

\def\NN{\mathbb{N}}

\begin{document}
\pagestyle{myheadings}
\markboth{\centerline{Benedek Schultz,~Jen\H{o} Szirmai }}
{Geodesic ball packings $\dots$,
}
\title
{Geodesic ball packings
generated by regular prism tilings in $\NIL$ geometry \footnote{AMS Classification 2000: 52C17, 52C22, 53A35, 51M20}}

\author{ B.~Schultz, ~ J. Szirmai\footnote{
 E-mail:
schultzb@math.bme.hu,
szirmai@math.bme.hu
}\\
\normalsize Budapest University of Technology and \\
\normalsize Economics Institute of Mathematics, \\
\normalsize Department of Geometry\\
}
%%%%%%%%%%%%%%%%%%%%%%%%%%%%%%%%%%%%%%%%%%%%
%\footnote{AMS Classification 2000: 52C17, 52C22, 53A35, 51M20}
%%%%%%%%%%%%%%%%%%%%%%%%%%%%%%%%%%%%%%%%%%%%

\maketitle
\begin{abstract}
In this paper we study the regular prism tilings and
construct ball packings by geodesic balls related to the above tilings in the projective model of $\NIL$ geometry.
Packings are generated by action of the discrete prism groups $\mathbf{pq2_{1}}$.
We prove that these groups are realized by prism tilings in $\NIL$ space if $(p,q)=(3,6), (4,4), (6,3)$ and determine
packing density formulae for geodesic ball packings generated by the above prism groups.
Moreover, studying these formulae we determine the conjectured maximal dense packing arrangements and their densities
and visualize them in the projective model of $\NIL$ geometry. We get a dense (conjectured locally densest) geodesic ball arrangement
related to the parameters $(p,q)=(6,3)$ where the kissing number of the packing is $14$, similarly to the densest lattice-like
$\NIL$ geodesic ball arrangement investigated by the second author in \cite{Sz07}.
\end{abstract}
%%%%%%%%%%%%%%%%%%%%%%%%%%%%%%%%%%%%%%%%%%%%

%%%%%%%%%%%%%%%%%%%%%%%%%%%%%%%%%%%%%%%%%%%
\newtheorem{theorem}{Theorem}[section]
\newtheorem{corollary}[theorem]{Corollary}
\newtheorem{lemma}[theorem]{Lemma}
\newtheorem{exmple}[theorem]{Example}
\newtheorem{definition}[theorem]{Definition}
\newtheorem{remark}[theorem]{Remark}
\newtheorem{conjecture}[theorem]{Conjecture}
\newtheorem{proposition}[theorem]{Proposition}
\newenvironment{example}{\begin{exmple}\normalfont}{\end{exmple}}
\newenvironment{acknowledgement}{Acknowledgement}

%%%%%%%%%%%%%%%%%%%%%%%%%%%%%%%%%%%%%%%%%%%%%%%%%%%%%%%%%%%%%%%%%%%%
%============================================================================%
%                             the main article                               %
%============================================================================%
%%%%%%%%%%%%%%%%%%%%%%%%%%%%%%%%%%%%%%%%%%%%%%%%%%%%%%%%%%%%%%%%%%%%

%%%%%%%%%%%%%%%%%%%%%%%%%%%%%%%%%%%%%%%%%%%%%%%%%%%%%%%%%%%%%%%%%%%%

%%%%%%%%%%%%%%%%%%%%%%%%%%%%%%%%%%%%%%%%%%%%%%%%%%%%%%%%%%%%%%%%%%%%%%%%%%%%%%
\section{Introduction and previous results}
In mathematics sphere packing problems concern the arrangements of non-over\-lapping equal spheres which fill a space.
Usually the space involved is the three-dimensional Euclidean space where the famous {\it Kepler conjecture} was proved by
T.~C.~Hales and S.~P.~Ferguson in \cite{HF06}.

However, ball (sphere) packing problems can be generalized to the other $3$-dimensional Thurston geometries.

In an $n$-dimensional space of constant curvature $\bE^n$, $\bH^n$, $\bS^n$ $(n\ge2)$ let $d_n(r)$ be the density
of $n+1$ spheres of radius $r$ mutually touching one another with respect to the simplex spanned by the centres of the spheres.
L.~Fejes T\'oth and H.~S.~M.~Coxeter conjectured that in an $n$-dimensional space of constant curvature the density of
packing spheres of radius $r$ can not exceed $d_n(r)$.
This conjecture has been proved by C.~Roger in the Euclidean space.
The 2-dimensional case has been solved by L.~Fejes T\'oth. In an $3$-dimensional space of constant curvature the
problem has been investigated by B\"or\"oczky and Florian in \cite{BF64} and it has been studied by K.~B\"or\"oczky
in \cite{B78} for $n$-dimensional space of constant curvature $(n\ge 4)$.

In \cite{KSz}, \cite{KSz1}, \cite{Sz12-1} and  \cite{Sz13-1} we have studied some new aspects of
the horoball and hyperball packings in $\mathbf{H}^n$ and we have observed that the ball, horoball and hyperball
packing problems are not settled yet in the $n$-dimensional $(n\ge3)$ hyperbolic space.

In \cite{Sz13-2} we generalized the above problem of finding the densest geodesic and translation ball (or sphere) packing to the other
$3$-dimensional homogeneous geometries (Thurston geometries)
$$
\SLR,~\NIL,~\SXR,~\HXR,~\SOL,
$$
and in the papers \cite{Sz10}, \cite{Sz07}, \cite{Sz11-1}, \cite{Sz12-2}, \cite{Sz13-2} we investigated
several interesting ball packing and covering problems in the above geometries. We described in $\SXR$ geometry (see \cite{Sz13-2})
a candidate of the densest geodesic and translation ball arrangement whose density is $\approx 0.8750$.

In this paper we consider the $\NIL$ geometry that can be derived from {W. Heisenberg}'s
famous real matrix group. This group provides a non-commutative translation group of an affine 3-space.
E. Moln\'ar proved in \cite{M97}, that the homogeneous 3-spaces
have a unified interpretation in the projective 3-sphere $\mathcal{PS}^3$.
In this work we will use this projective model of the $\NIL$ geometry.

In \cite{Sz07} we investigated the geodesic balls of the $\NIL$ space and computed their volume,
introduced the notion of the $\NIL$ lattice, $\NIL$ parallelepiped and the density of the lattice-like ball packing.
Moreover, we determined the densest lattice-like geodesic ball packing. The density of this densest packing is
$\approx 0.7809$, may be surprising enough
in comparison with the Euclidean result $\frac{\pi}{\sqrt{18}} \approx 0.74048$. The kissing number of the balls
in this packing is $14$.

In \cite{Sz12-2} we considered the analogue question for translation balls. The notions of translation curve and
translation ball were introduced by initiative of {E. Moln\'ar} (see \cite{MSz}, \cite{Sz10}).
We have studied the translation balls of $\NIL$ space and computed their volume. Moreover, we have proved that the density
of the optimal lattice-like translation ball
packing for every natural lattice parameter $ 1 \le k \in \mathbf{N}$ is in interval $(0.7808, 0.7889)$
and if  $r \in (0,r_{d}] \ (r_{d} \approx 0.7456)$
then the optimal density is $\delta^{opt} \approx  0.7808$. Meanwhile we can apply a nice general estimate of
{L.~Fejes T\'oth} \cite{FTL64} in our proof.
The kissing number of the lattice-like ball packings is less than or equal to 14 and the optimal ball packing is realizable
in case of equality. We formulated a conjecture for $\delta^{opt}$, where the density of the conjectural densest packing
is $\delta^{opt} \approx 0.7808$ for lattice parameter $k=1$, larger than the Euclidean one
($\frac{\pi}{\sqrt{18}} \approx 0.7405$), but
less than the density of the densest lattice-like {\it geodesic} ball packing in $\NIL$ space known till now \cite{Sz07}.
The kissing number of the translation balls in that packing is 14 as well.

In \cite{Sz13-2} we studied one type of {\it lattice coverings} in the $\NIL$ space. We introduced
the notion of the density of considered coverings and
gave upper and lower estimation to the density of the
lowest lattice-like geodesic ball covering. Moreover we formulate a conjecture for the ball
arrangement of the least dense lattice-like geodesic
ball covering and give its covering density $\Delta\approx 1.4290$.

In papers \cite{Sz12-2},\cite{MSz} there are further results related to discrete problems of $\NIL$ geometry.

{\it In this paper we study the regular prism tilings and construct ball packings by geodesic balls related to the prism tilings
in the projective model of $\NIL$ geometry where the packings are generated by action of the discrete
prism groups  $\mathbf{pq2_{1}}$. We obtain density formulae for calculations for geodesic ball packings.
Analyzing these density functions we obtain a conjecture for optimal geodesic ball packing configurations
and determine their densities related to the above prismatic tessellations.

The results are summarized in Theorems 3.4, 4.3 and Conjecture \ref{conj:ball}, and the optimal ball configurations are visualized
in our model with Figures 2, 3.}
%
%%%%%%%%%%%%%%%%%%%%%%%%%%%%%%%%%%%%%%%%%%%%%%%%%%%%%%%%%%%%%%%%%%%%%%%%%%%%%%%%%%%%%%%%%%%%%
\section{Basic notions of the $\NIL$ geometry}
In this Section we summarize the significant notions and denotations of the $\NIL$ geometry (see \cite{M97}, \cite{Sz07}).

The $\NIL$ geometry is a homogeneous 3-space derived from the famous real matrix group $\mathbf{L(R)}$ discovered by {W.~Heisenberg}.
The Lie Theory with the method of the projective geometry makes possible to investigate and to describe this topic.

The left (row-column) multiplication of Heisenberg matrices
     \begin{equation}
     \begin{gathered}
     \begin{pmatrix}
         1&x&z \\
         0&1&y \\
         0&0&1 \\
       \end{pmatrix}
       \begin{pmatrix}
         1&a&c \\
         0&1&b \\
         0&0&1 \\
       \end{pmatrix}
       =\begin{pmatrix}
         1&a+x&c+xb+z \\
         0&1&b+y \\
         0&0&1 \\
       \end{pmatrix}
      \end{gathered} \tag{2.1}
     \end{equation}
defines "translations" $\mathbf{L}(\mathbf{R})= \{(x,y,z): x,~y,~z\in \mathbf{R} \}$
on the points of the space $\NIL= \{(a,b,c):a,~b,~c \in \mathbf{R}\}$.
These translations are not commutative in general. The matrices $\mathbf{K}(z) \vartriangleleft \mathbf{L}$ of the form
     \begin{equation}
     \begin{gathered}
       \mathbf{K}(z) \ni
       \begin{pmatrix}
         1&0&z \\
         0&1&0 \\
         0&0&1 \\
       \end{pmatrix}
       \mapsto (0,0,z)
      \end{gathered}\tag{2.2}
     \end{equation}
constitute the one parametric centre, i.e. each of its elements commutes with all elements of $\mathbf{L}$.
The elements of $\mathbf{K}$ are called {\it fibre translations}. $\NIL$ geometry of the Heisenberg group can be projectively
(affinely) interpreted by the "right translations"
on points as the matrix formula
     \begin{equation}
     \begin{gathered}
       (1;a,b,c) \to (1;a,b,c)
       \begin{pmatrix}
         1&x&y&z \\
         0&1&0&0 \\
         0&0&1&x \\
         0&0&0&1 \\
       \end{pmatrix}
       =(1;x+a,y+b,z+bx+c)
      \end{gathered} \tag{2.3}
     \end{equation}
shows, according to (1.1). Here we consider $\mathbf{L}$ as projective collineation
group with right actions in homogeneous coordinates.
We will use the Cartesian homogeneous coordinate simplex $E_0(\be_0)$,$E_1^{\infty}(\be_1)$,$E_2^{\infty}(\be_2)$,
$E_3^{\infty}(\be_3), \ (\{\be_i\}\subset \bV^4$ \ $\text{with the unit point}$ $E(\be = \be_0 + \be_1 + \be_2 + \be_3 ))$
which is distinguished by an origin $E_0$ and by the ideal points of coordinate axes, respectively.
Moreover, $\by=c\bx$ with $0<c\in \mathbb{R}$ (or $c\in\mathbb{R}\setminus\{0\})$
defines a point $(\bx)=(\by)$ of the projective 3-sphere $\cP \cS^3$ (or that of the projective space $\cP^3$ where opposite rays
$(\bx)$ and $(-\bx)$ are identified).
The dual system $\{(\Be^i)\}, \ (\{\be^i\}\subset \BV_4)$ describes the simplex planes, especially the plane at infinity
$(\Be^0)=E_1^{\infty}E_2^{\infty}E_3^{\infty}$, and generally, $\Bv=\Bu\frac{1}{c}$ defines a plane
$(\Bu)=(\Bv)$ of $\cP \cS^3$
(or that of $\cP^3$). Thus $0=\bx\Bu=\by\Bv$ defines the incidence
of point $(\bx)=(\by)$ and plane
$(\Bu)=(\Bv)$, as $(\bx) \text{I} (\Bu)$ also denotes it.
Thus {\bf Nil} can be visualized in the affine 3-space $\bA^3$
(so in $\bE^3$) as well.

The translation group $\mathbf{L}$ defined by formula (2.3) can be extended to a larger group $\mathbf{G}$ of collineations,
preserving the fibering, that will be equivalent to the (orientation preserving) isometry group of $\NIL$.
In \cite{M06} E.~Moln\'ar has shown that
a rotation trough angle $\omega$
about the $z$-axis at the origin, as isometry of $\NIL$, keeping invariant the Riemann
metric everywhere, will be a quadratic mapping in $x,y$ to $z$-image $\overline{z}$ as follows:
     \begin{equation}
     \begin{gathered}
       \br(O,\omega):(1;x,y,z) \to (1;\overline{x},\overline{y},\overline{z}); \\
       \overline{x}=x\cos{\omega}-y\sin{\omega}, \ \ \overline{y}=x\sin{\omega}+y\cos{\omega}, \\
       \overline{z}=z-\frac{1}{2}xy+\frac{1}{4}(x^2-y^2)\sin{2\omega}+\frac{1}{2}xy\cos{2\omega}.
      \end{gathered} \tag{2.4}
     \end{equation}
This rotation formula, however, is conjugate by the quadratic mapping $\mathcal{M}$
     \begin{equation}
     \begin{gathered}
       x \to x'=x, \ \ y \to y'=y, \ \ z \to z'=z-\frac{1}{2}xy  \ \ \text{to} \\
       (1;x',y',z') \to (1;x',y',z')
       \begin{pmatrix}
         1&0&0&0 \\
         0&\cos{\omega}&\sin{\omega}&0 \\
         0&-\sin{\omega}&\cos{\omega}&0 \\
         0&0&0&1 \\
       \end{pmatrix}
       =(1;x",y",z"), \\
       \text{with} \ \ x" \to \overline{x}=x", \ \ y" \to \overline{y}=y", \ \ z" \to \overline{z}=z"+\frac{1}{2}x"y",
      \end{gathered} \tag{2.5}
     \end{equation}
i.e. to the linear rotation formula. This quadratic conjugacy modifies the $\NIL$ translations in (2.3), as well.
We shall use the following important classification theorem.
\begin{theorem}[E.~Moln\'ar \cite{M06}]
\begin{enumerate}
	\item Any group of $\NIL$ isometries, containing a 3-dimensional translation lattice,
is conjugate by the quadratic mapping in (2.5) to an affine group of the affine (or Euclidean) space $\bA^3=\EUC$
whose projection onto the (x,y) plane is an isometry group of $\bE^2$. Such an affine group preserves a plane
$\to$ point polarity of signature $(0,0,\pm 0,+)$.
	\item Of course, the involutive line reflection about the $y$ axis
     \begin{equation}
     \begin{gathered}
       (1;x,y,z) \to (1;-x,y,-z),
      \end{gathered} \notag
     \end{equation}
preserving the Riemann metric, and its conjugates by the above isometries in {\bf 1} (those of the identity component)
are also {$\NIL$}-isometries. There does not exist orientation reversing $\NIL$-isometry.
\end{enumerate}
\end{theorem}
\begin{remark}
We obtain from the above described projective model a new model of $\NIL$ geometry derived by the quadratic mapping $\mathcal{M}$.
This is the {\it linearized model of $\NIL$ space} (see \cite{Br15}).
\end{remark}
\subsection{Geodesic curves and spheres}
The geodesic curves of the $\NIL$ geometry are generally defined as having locally minimal arc length between their any two (near enough) points.
The equation systems of the parametrized geodesic curves $g(x(t),y(t),z(t))$  in our model can be determined by the
general theory of Riemann geometry.
We can assume, that the starting point of a geodesic curve is the origin because we can transform a curve into an
arbitrary starting point by translation (2.1);
\begin{equation}
\begin{gathered}
        x(0)=y(0)=z(0)=0; \ \ \dot{x}(0)=c \cos{\alpha}, \ \dot{y}(0)=c \sin{\alpha}, \\ \dot{z}(0)=w; \ - \pi \leq \alpha \leq \pi. \notag
\end{gathered}
\end{equation}
The arc length parameter $s$ is introduced by
\begin{equation}
 s=\sqrt{c^2+w^2} \cdot t, \ \text{where} \ w=\sin{\theta}, \ c=\cos{\theta}, \ -\frac{\pi}{2}\le \theta \le \frac{\pi}{2}, \notag
\end{equation}
i.e. unit velocity can be assumed.
\begin{remark}
Thus we have harmonized the scales along the coordinate axes.
\end{remark}
The equation systems of a helix-like geodesic curves $g(x(t),y(t),z(t))$ if $0<|w| <1 $:
\begin{equation}
\begin{gathered}
x(t)=\frac{2c}{w} \sin{\frac{wt}{2}}\cos\Big( \frac{wt}{2}+\alpha \Big),\ \
y(t)=\frac{2c}{w} \sin{\frac{wt}{2}}\sin\Big( \frac{wt}{2}+\alpha \Big), \notag \\
z(t)=wt\cdot \Big\{1+\frac{c^2}{2w^2} \Big[ \Big(1-\frac{\sin(2wt+2\alpha)-\sin{2\alpha}}{2wt}\Big)+ \\
+\Big(1-\frac{\sin(2wt)}{wt}\Big)-\Big(1-\frac{\sin(wt+2\alpha)-\sin{2\alpha}}{2wt}\Big)\Big]\Big\} = \\
=wt\cdot \Big\{1+\frac{c^2}{2w^2} \Big[ \Big(1-\frac{\sin(wt)}{wt}\Big)
+\Big(\frac{1-\cos(2wt)}{wt}\Big) \sin(wt+2\alpha)\Big]\Big\}. \tag{2.6}
\end{gathered}
\end{equation}
In the cases $w=0$ the geodesic curve is the following:
\begin{equation}
x(t)=c\cdot t \cos{\alpha}, \ \ y(t)=c\cdot t \sin{\alpha}, \ \ z(t)=\frac{1}{2} ~ c^2 \cdot t^2 \cos{\alpha} \sin{\alpha}. \tag{2.7}
\end{equation}
The cases $|w|=1$ are trivial: $(x,y)=(0,0), \ z=w \cdot t$.
\begin{definition}
The distance $d(P_1,P_2)$ between the points $P_1$ and $P_2$ is defined by the arc length of geodesic curve
from $P_1$ to $P_2$.
\end{definition}
In our work \cite{Sz07} we introduced the following definitions:
\begin{definition}
The geodesic sphere of radius $R$ with centre at the point $P_1$ is defined as the set of all points
$P_2$ in the space with the condition $d(P_1,P_2)=R$. Moreover, we require that the geodesic sphere is a simply connected
surface without self-intersection
in the $\NIL$ space.
\begin{remark}
We shall see that this last condition depends on radius $R$.
\end{remark}
\end{definition}
\begin{definition}
The body of the geodesic sphere of centre $P_1$ and of radius $R$ in the $\NIL$ space is called geodesic ball, denoted by $B_{P_1}(R)$,
i.e. $Q \in B_{P_1}(R)$ iff $0 \leq d(P_1,Q) \leq R$.
\end{definition}
\begin{remark}
Henceforth, typically we choose the origin as centre of the sphere and its ball, by the homogeneity of
$\NIL$.
\end{remark}
We have denoted by $B(S)$ the body of the $\NIL$ sphere $S$,
furthermore we have denoted their volumes by $Vol(B(S))$.

In \cite{Sz07} we have proved the the following theorem:
\begin{theorem}
The geodesic sphere and ball of radius $R$ exists in the $\NIL$ space if and only if $R \in [0,2\pi].$
\end{theorem}
We obtain the volume of the geodesic ball of radius $R$ by the following integral (see 2.8):
\begin{equation}
\begin{gathered}
Vol(B(S))=2 \pi \int_0^{\frac{\pi}{2}}X^2 ~ \frac{\mathrm{d}~Z}{\mathrm{d}~\theta} ~ \mathrm{d}~{\theta}= \\ =2 \pi \int_0^{\frac{\pi}{2}}
\Big (\frac{2\cos{\theta}}{\sin{\theta}} \sin{\frac{(R \sin{\theta})}{2}}\Big)^2 \cdot \Big(-\frac{1}{2} \frac{R \cos^3{\theta}}{\sin^2{\theta}}+
\frac{\cos{\theta}\sin{(R \sin{\theta})}}{\sin{\theta}}+ \\ +\frac{\cos^3{\theta}\sin{(R \sin{\theta})}}{\sin^3{\theta}}-
\frac{1}{2} \frac{R \cos^3{\theta}\cos{(R \sin{\theta})}}{\sin^2{\theta}}\Big){\mathrm{d}~\theta}. \tag{2.8}
\end{gathered}
\end{equation}
The parametric equation system of the geodesic sphere $S(R)$ in
our model (see \cite{Sz07}):
\begin{equation}
\begin{gathered}
x(R,\theta,\phi)=\frac{2c}{w} \sin{\frac{wR}{2}}\cdot \cos{\phi}= \frac{2\cos{\theta}}{\sin{\theta}}
\sin{\frac{R \sin{\theta}}{2}}\cdot \cos{\phi}, \notag \\
y(R,\theta,\phi)=\frac{2c}{w} \sin{\frac{wR}{2}}\cdot \sin{\phi}= \frac{2\cos{\theta}}{\sin{\theta}}
\sin{\frac{R \sin{\theta}}{2}}\cdot \sin{\phi}, \notag \\
z(R,\theta,\phi)=wR+\frac{c^2R}{2w}-\frac{c^2}{2w^2}\sin{wR} +\frac{1}{4}\Big( \frac{2c}{w}
\sin{\frac{wR}{2}}\Big)^2 \sin{2\phi}= \\
=R\sin{\theta}+\frac{R\cos^2{\theta}}{2\sin{\theta}} - \frac{\cos^2{\theta}}{2\sin^2{\theta}}\sin(R\sin{\theta})+
\frac{1}{4}\Big( \frac{2\cos{\theta}}{\sin{\theta}} \sin{R \frac{\sin{\theta}}{2}}\Big)^2 \sin{2\phi} \notag \\
\end{gathered}
\end{equation}
 \begin{equation}
\begin{gathered}
  -\pi<\phi\leqq \pi, \ \ -\frac{\pi}{2}\leqq \theta \leqq \frac{\pi}{2} \ \text{and} \ \theta \ne 0. \\
   \text{if} \ \theta=0 \ \text{then} \ x(R,0,\phi)=R \cos{\phi}, \ \ y(R,0,\phi)=R \sin{\phi}, \\
   \ z(R,0,\phi)=\frac{1}{2} ~ R^2 \cos{\phi} \sin{\phi}. \tag{2.9}
\end{gathered}
\end{equation}
We have obtained by the derivatives of these parametrically represented functions (by intensive and careful computations with {\it Maple}
through the second fundamental form) the following theorem (see \cite{Sz07}):
\begin{theorem}
The geodesic $\NIL$ ball $B(S(R))$ is convex in affine-Euclidean sense in our model if and only if $R \in [0,\frac{\pi}{2}]$.
\end{theorem}
%

%%%%%%%%%%%%%%%%%%%%%%%%%%%%%%%%%%%%%%%%%%%%%%%%%%%%%%%%%%%%%%%%%%%%%%%%%%%%%%%%%%%%
\section{$\NIL$ prisms and prism tilings}
The prisms and prism-like tilings have been thoroughly investigated in $\SXR, \HXR$ and
$\SLR$ spaces in papers \cite{PSchSz12}, \cite{SchSz15}, \cite{Sz14-2}.
Here we consider the analogous problem in $\NIL$ space.
We will use the in 2. section described projective model of $\NIL$ geometry.
In the following the plane of $x,~y$ axis are called {\it base plane}
of the model and if we say {\it plane} then it is a plane in Euclidean sense.
\begin{definition}
Let $\cP^i$ be an infinite solid bounded by planes, that are determined by fibre-lines passing through the points
of a $p$-gon ($p\ge3$, integer parameter) $\cP^b$ lying in the base-plane. The images of $\cP^i$
by $\NIL$ isometries are called infinite $p$-sided prisms.

The common part of $\cP^i$ with the base plane is defined as the \textit{base figure} $\cP^b$ of the prism.
\end{definition}
Let $\cF$ be the $\cM^{-1}$ image of the base plane in the $\NIL$-space (see Remark 2.2)
and let $\tau$ be a fibre translation (2.2).

\begin{definition}
Let $\cP^i$ be an infinite $p$-sided prism, that is trimmed by the surface $\cF$ and its translated copy $\cF^\tau$.
The parts of $\cF$ and $\cF^\tau$ inside the infinite prism are called cover faces and are
denoted by $C_{\cF}$ and $C_{\cF^\tau}$.

The $p$-sided bounded prism is the part of $\cP^i$ between the cover faces $C_{\cF}$ and $C_{\cF^\tau}$.
\end{definition}
\begin{definition}
A bounded or infinite $p$-sided prism is said to be \textit{regular} if its side surfaces are congruent to each other under
$\NIL$ rotations with angle $\frac{2\pi}{p}$ (see (2.4) and (2.5)) about the central fibre line of the prism.
\end{definition}
\subsection{Regular bounded prism tilings}
In this section we will investigate the existence of regular bounded prism tilings $\cT_p(q)$ of $\NIL$ space.
In this case the prism tiles are regular bounded prisms having $p$-gonal base figures $(p \ge 3)$.
The prism itself is a \textit{topological polyhedron} with $2p$ vertices, and having at every vertex one
$p$-gonal cover face and two quadrangle side faces (traced by fibre lines). We are looking such prism tilings of
$\NIL$ space where at each side edge of the prism
(which are fibre lines going through vertices of the base figure) meet $q$ prisms regularly,
by $\NIL$ rotations with angle $\frac{2\pi}{q}$ ($q\ge 3$, integer parameter).

We shall see in Theorem 3.4 that the regular prism tiling $\cT_p(q)$ exists for some parameters $(p,q)$. Let $\cP_p(q)$ one of its tiles with
with vertices $A_1A_2 \dots A_p$ $B_1B_2 \dots B_p$. We may assume that $A_1$ lies on the
$x$-axis. It is clear that the side curves $c_{A_iA_{i+1}}$ $(i=1\dots p, ~ A_{p+1} \equiv A_1)$
are derived from each other by $\frac{2\pi}{p}$ rotation about the $x$ axis.
The corresponding vertices $B_1B_2 \dots B_p$ are generated by a fibre translation $\tau$
with a positive real parameter. The cover faces $A_1,\dots, A_p$, $B_1,\dots, B_p$ and the side surfaces form a
$p$-sided regular prism $\cP_p(q)$ in $\NIL$.
$\cT_p(q)$ will be generated by its rotational isometry group $\Gamma_p(q)=\mathbf{pq2_1}$ $\subset Isom(\NIL)$
(if these tiling there exist see Theorem 3.4) which is given by its fundamental domain
$\mathcal{F}_p(q)=A_1A_2O A_1^{\bs} A_2^{\bs} O^{\bs}$, $A_1^{\bs}=B_p$. Here,
$A_2^{\bs}=B_1$, $O^{\bs}=O^{\tau}$, and $\mathcal{F}_p(q)$ is
a piece-wise linear topological polyhedron. The group presentation can be determined
by a standard procedure \cite{M92}, called Poincar\'e
algorithm. The generators will pair the bent (piecewise linear) faces of $\mathcal{F}$:
\begin{equation}
\begin{gathered}
\ba~:~OA_1B_pO^{\bs} (O) \rightarrow OA_2B_1O^{\bs} (O), \\
\bb~:~A_1A_2B_1 (A_1) \rightarrow A_1B_pB_1 (A_1),~
\bs~:~OA_1A_2 (O) \rightarrow O^{\bs}B_pB_1(O^{\bs})
\end{gathered} \notag
\end{equation}
mapping $\mathcal{F}_p(q)$ onto its neighbours ${\mathcal{F}_p(q)}^{\ba}$, ${\mathcal{F}_p(q)}^{\bb}$,
${\mathcal{F}_p(q)}^{\bs}$, respectively.
E.g. for the face $a^{-1}$ a point $A$ (relative freely, e.g. in the segment $OB_p$) is taken.
Then the union of triangles $AOO^{\bs}$, $AO^{\bs}B_p$, $AB_pA_1$, $OA_1O$ will be the face $a^{-1}$.

Then the $\ba$-image $A^{\ba}$ is taken in $OB_1$ for the face $a=A^{\ba}OO^{\bs} \cup
A^{\ba}O^{\bs}B_1 \cup A^{\ba}B_1A_2 \cup A^{\ba}A_2O$, as usual. The relations are induced by the edge equivalence classes $\{OO'\}$; $\{A_1B_1\}$; $\{OA_1,$ $OA_2,O'B_1,O'B_p\}$;
$\{A_1A_2,A_1B_p,A_2B_1,B_pB_1 \}$. So we get the group
\begin{equation}
\mathbf{pq2_1}=\{\ba,\bb:\ba^p=\bb^q=\ba\bb\ba\bb\ba^{-1}\bb^{-1}\ba^{-1}\bb^{-1}=1\}, \tag{3.1}
\end{equation}
where $\ba$ is a $p$ rotation about the fibre line through the origin (the $z$-axis),
$\bb$ is a $q$ rotation about a side fibre line of $\cP_p(q)$ (through a vertex of its base figure).
Notice, that $\bb\ba\bb$ is a screw motion, and thus $\tau:=\ba\bb\ba\bb=\bb\ba\bb\ba$ is
the fibre translation connecting the cover faces.

Our first question is the following: For which $3 \le p,q \in\NN$ is $\Gamma_p(q)=\mathbf{pq2_1}\subset\mathit{Isom}(\NIL)$?

The following Theorem answers it:
\begin{theorem}\label{thm:geodb}
 In $\NIL$ there exist $3$ regular $p$-gonal non-face-to-face prism tilings $\cT_p(q)$ with $\NIL$ isometry group
 $\Gamma_p(q)=\mathbf{pq2_1}$ for integer parameters $p,q \ge 3$:

the regular triangular prism tiling with $(p,q)=(3,6)$,

the regular square prism tiling with $(p,q)=(4,4)$,

the regular hexagonal prism tiling with $(p,q)=(6,3)$,

\noindent and each group $\Gamma_p(q)$ has a free parameter $x_p(q) \in \mathbb{R}^+$.
\end{theorem}
\begin{figure}[htbp]\label{pic:prism1}
\centering
\includegraphics[height=6cm]{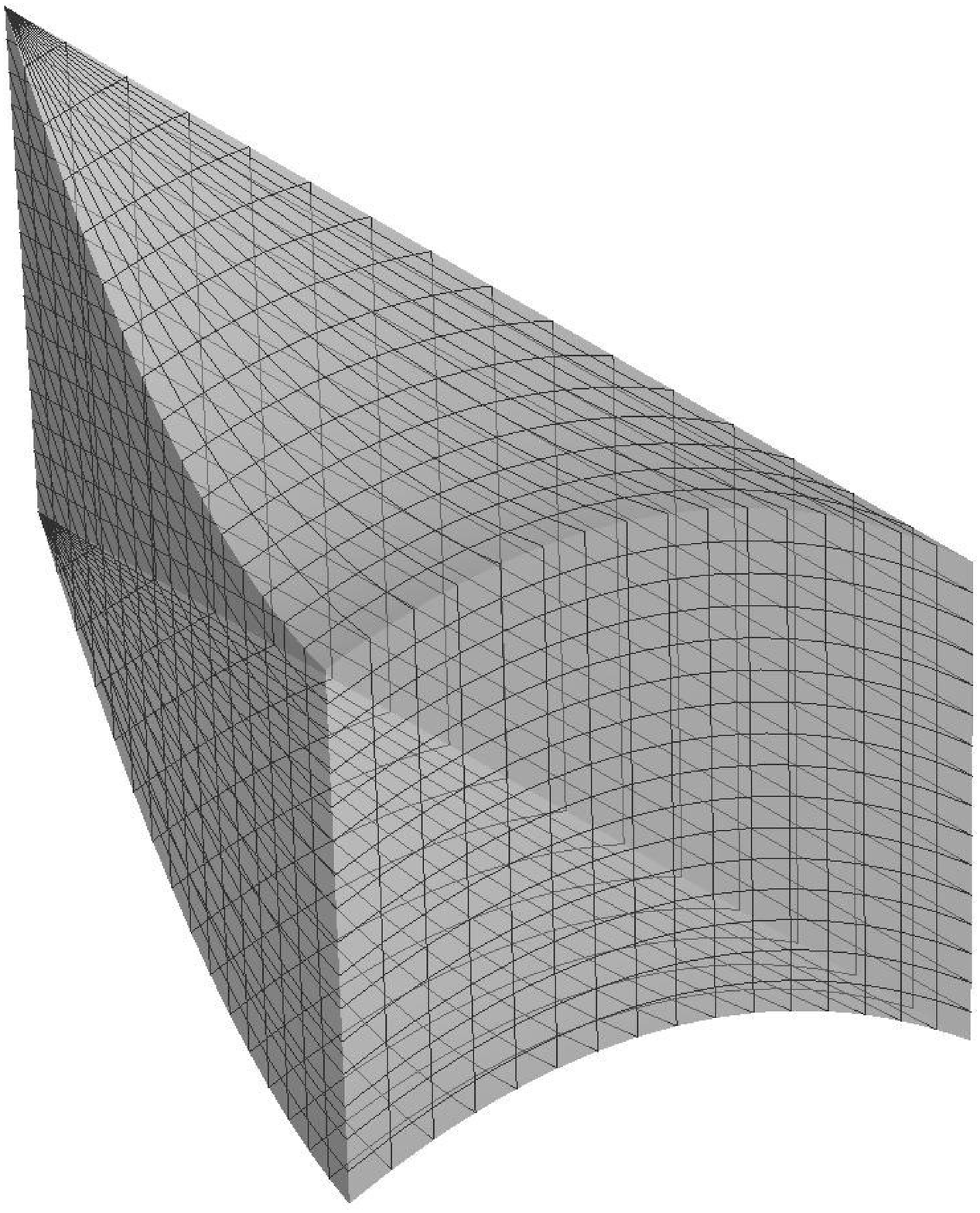}\includegraphics[height=6cm]{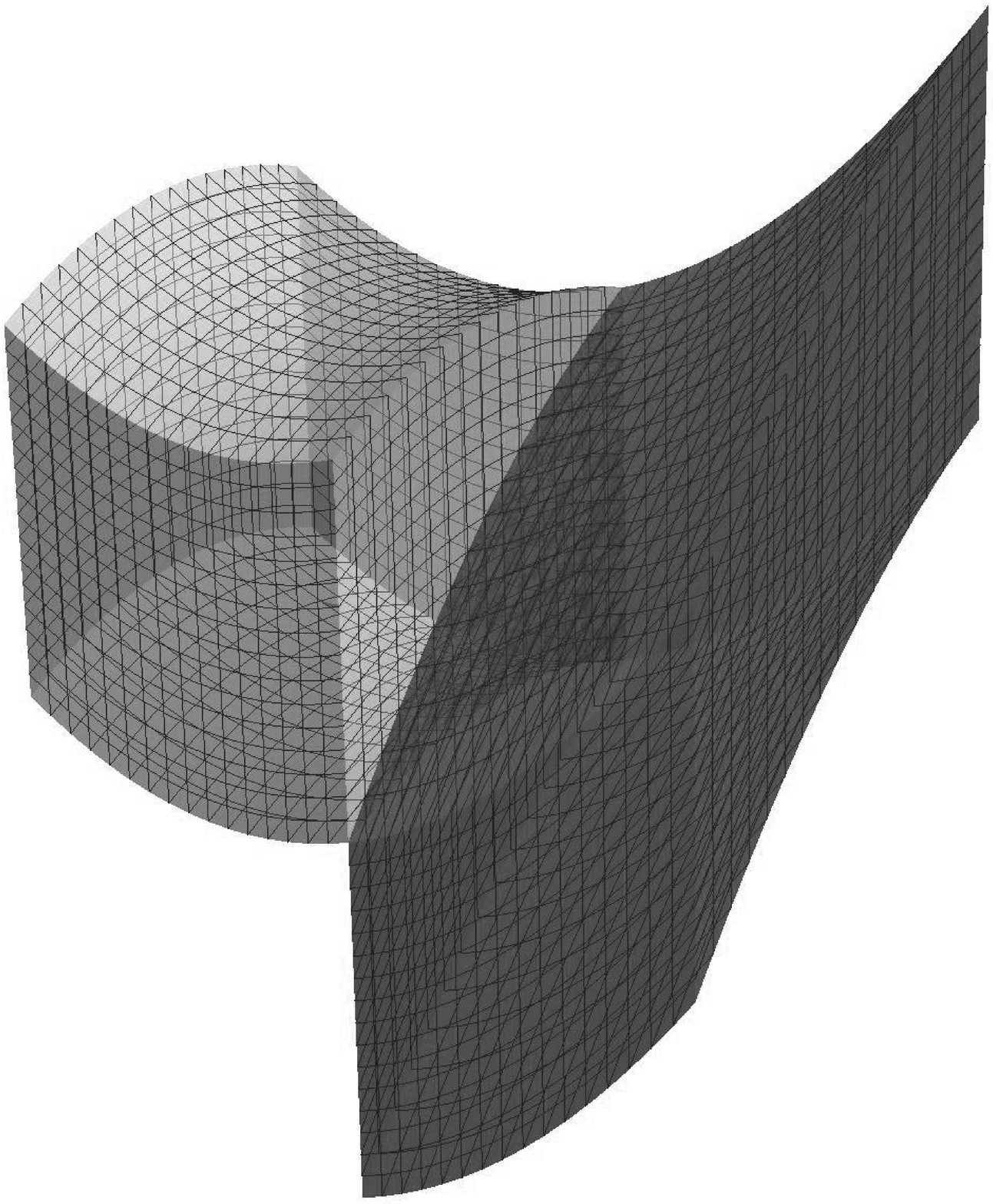}
\caption{The regular triangular and rectangular prisms in $\NIL$-space, with parameters $x_3(6) = \frac{2}{\sqrt{3}}$ and
$x_4(4) = \sqrt{2}$ respectively}
\label{pic:36}
\end{figure}

{\bf Proof}:

Let $A_1=(1;x_p(q),0,0)$ be a "bottom" vertex of the regular bounded prism. Then the other vertices of the bottom cover face can be generated
by the $\NIL$ rotation formula (see (2.4), (2.5)):
\begin{gather*}
 A_2=A_1^\ba=\Big(1; x_p(q)\cos\Big(\frac{2\pi}{p}\Big), x_p(q) \sin\Big(\frac{2\pi}{p}\Big), \frac{1}{4}{x_p(q)}^2\sin\Big(\frac{4\pi}{p}\Big)\Big),\\
 A_3=A_2^\ba=A_1^{\ba^2}=\Big(1; x_p(q)\cos\Big(\frac{4\pi}{p}\Big), x_p(q)\sin\Big(\frac{4\pi}{p}\Big), \frac{1}{4}x_p(q)^2\sin\Big(\frac{8\pi}{p}\Big)
 \Big),\\ \dots\\
 A_p=A_{p-1}^\ba=A_{p-1}^{\ba^{p-1}}.
\end{gather*}

Then the condition for the existence of the tiling is the following:
\begin{equation}
 A_3=A_2^\ba\equiv A_1^{\bb^{-1}},\tag{3.2}
\end{equation}
where $\ba$ is a $p$ rotation about the fibre line through the origin and
$\bb$ is a $q$ rotation about the side fibre line of $\cP_p(q)$ through the vertex $A_2$.
$\equiv$ means that the corresponding points lie on the same fibre lines.
\begin{equation*}
 1=\cos^2\Big(\frac{\pi}{p}\Big)+\cos^2\Big(\frac{\pi}{q}\Big),
\end{equation*}
where $p$ and $q$ are positive integers. This equation only has the following integer solutions:
\begin{equation*}
 (p,q)=(4,4),(3,6)\text{ or }(6,3).
\end{equation*}
We obtain from the above computations, that the existence of the above regular prism tilings is independent from the parameter
$x_p(q) \in \mathbb{R}^+$, so we have proven the Theorem. ~ $\square$

Remembering that $\tau=\ba\bb\ba\bb$ is the "vertical" translation of the group, we can also compute the height of the regular bounded prism corresponding to the group tiling, since:
$O^{\ba\bb\ba\bb}=O^\tau,$
where $O$ is the origin.
Using this, we can also give a metric representation of the group,
allowing the visualization of the corresponding prism and prism tiling (see Fig.~1. and Fig.~2.).
\section{The optimal geodesic ball packings under group $\mathbf{pq2_1}$}
The sphere packing problem deals with the arrangements of non-overlapping equal spheres,
or balls, which fill the space. While the usual problem is in the $n$-dimensional Euclidean-space $(n\ge 2$),
it can be generalized to the other $3$-dimen\-si\-o\-nal Thurston spaces (see \cite{Sz14-1}).
In this paper we investigate the optimal ball packings of $\NIL$ generated by the above described $\mathbf{pq2_1}$ group.

Let $\cT_p(x_p(q),q)$ (where $(p,q)=(3,6),(4,4)$ or $(6,3)$ as stated above and $x_p(q)\in \mathbb{R}^+$) be a regular prism tiling,
and let $\cP_p(x_p(q),q)$ be one of its tiles that is centered at the origin, with a base face given by the vertices
$A_1,A_2,\dots,A_p$. The corresponding vertices $B_1,B_2,\dots,B_p$ of the prism are generated by fibre translations $\tau=\ba\bb\ba\bb$.

We can assume by symmetry, that the optimal geodesic ball is centered at the origin.
The volume of a geodesic ball with radius $R$ can be determined by the formula (2.8).

We study only one case of the multiply transitive geodesic ball packings where the fundamental domains
of the $\NIL$ space groups $\mathbf{pq2_1}$ are not prisms.
Let the fundamental domains be derived by the Dirichlet~---~Voronoi cells (D-V cells) where their centers are images of the origin.
The volume of the $p$-times fundamental domain and of the D-V cell is the same, respectively,
as in the prism case (for any above $(p, q, x_p(q))$ fixed). It is easy to see by the formulas (2.5),
using the quadratic mapping $\cM$, that the volume of the Dirichlet~---~Voronoi cell (or the coresponding prism) is
\begin{equation}
Vol(\cP_p(x_p(q),q))=\frac{p}{2} {x_p^2(q)} \sin\Big( \frac{2\pi}{p} \Big)d(OO^\tau). \tag{4.1}
\end{equation}
These locally densest geodesic ball packings can be determined for all possible fixed integer parameters $p,q, x_p(q)$.
The optimal radius $R_{opt}(x_p(q))$ is
\begin{equation}
R_{opt}(x_p(q),p,q) = \min \Big\{ d{(OA_1)}, \ \frac{d(OO^\tau)}{2}, \ \frac{d(O,O^{{\mathbf a} {\mathbf b} })}{2} \Big\} , \tag{4.2}
\end{equation}
where $d$ is the geodesic distance function of $\NIL$ geometry (see Definition 2.4).

Since the congruent images of $\cP_p(x_p(q),q)$ under the discrete group $\mathbf{pq2_1}$
cover the $\NIL$ space, therefore for the density of the ball packing
it is sufficient to relate the volume of the ball to the volume of the prism:
\begin{definition}
The maximal density $\delta_p(x_p(q),q)$ of the above multiply transitive ball packing for given parameters $(p, q, x_p(q))$ ($(p,q)=(3,6),(4,4),(6,3)$
and $x_p(q)\in \mathbb{R}^+$):
\begin{equation}
 \delta_p(x_p(q),q)=\frac{Vol(B(R_{opt}))}{Vol(\cP_p(x_p(q),q))}=\frac{Vol(B(R_{opt}))}{\frac{p}{2} {x_p^2(q)}
 \sin\Big( \frac{2\pi}{p} \Big)d(OO^\tau)}. \tag{4.3}
\end{equation}
\end{definition}
For every $p,q,x_p(q)$ parameters the locally densest geodesic ball packing can be determined.

If we fixed the parameters $p$ and $q$ then the distance function $d(x_p(q))$ is a continuous fuction. Therefore it is easy prove
the following Lemma:
\begin{lemma}
In $\NIL$-space for the rotation group $\Gamma_p(q)=\mathbf{pq2_1}$ there always exist $x_p(q)\in \mathbb{R}^+$
for given parameters $(p,q)=(4,4),(3,6),(6,3)$ where
\begin{equation}\label{eq:optcond}
 d(O, O^{\ba\bb})=d(O,O^{\ba\bb\ba\bb})=d(O,O^{\tau}).\tag{4.4}
\end{equation}
\end{lemma}
\begin{figure}[htbp]\label{pic:prism3}
\centering
\includegraphics[height=4.5cm]{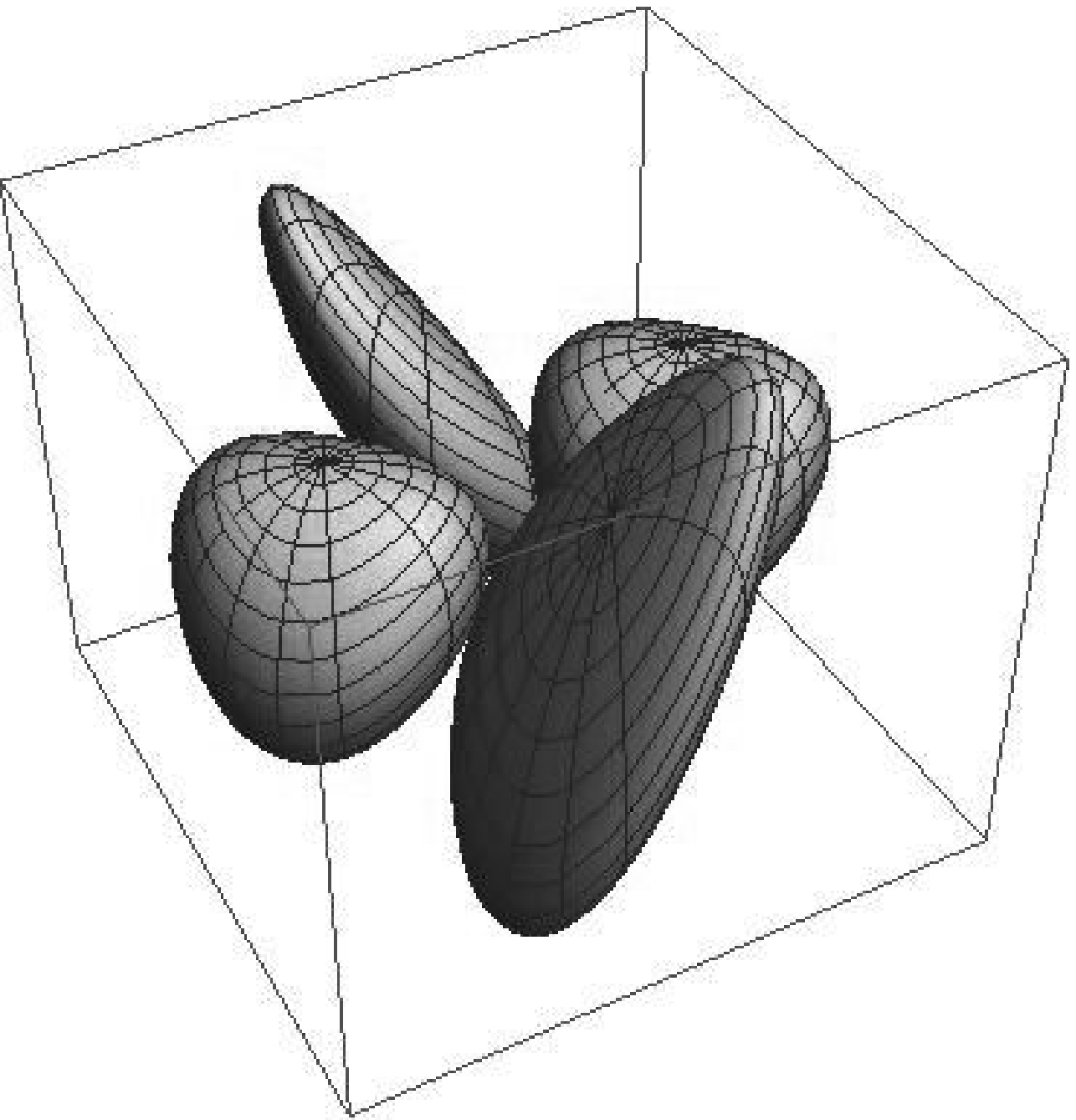}\includegraphics[height=7.5cm]{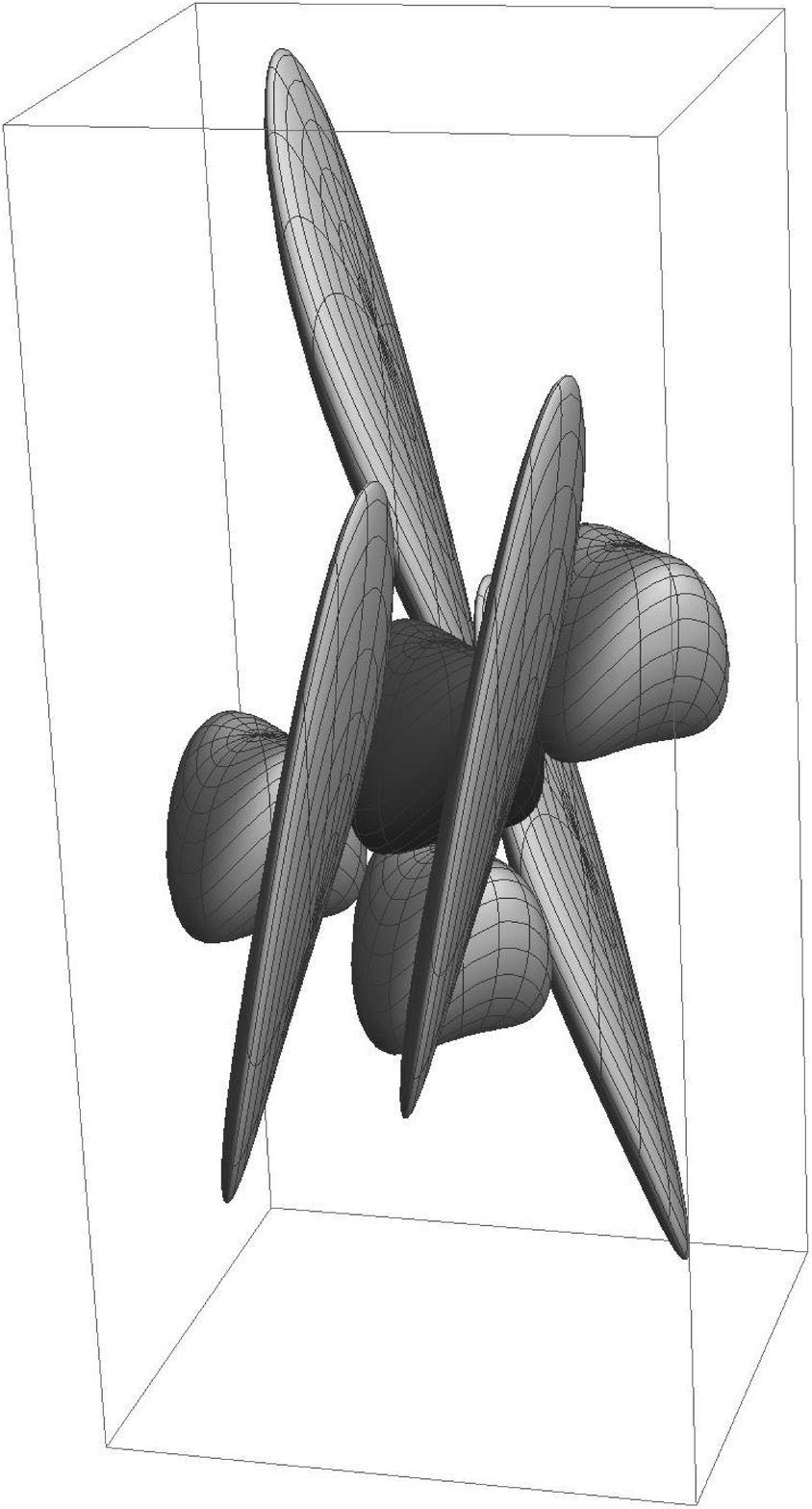}
\caption{Some balls of the optimal ball arrangements for the square and hexagonal tilings with parameters $(p,q)=(4,4)$
and $(p,q)=(3,6)$.}
\label{pic:optballs}
\end{figure}
{\it The system of equations (4.4) in Lemma 4.2 can be solved by numerical methods and the corresponding
ball arrangements are denoted by $\cB_p(q)$}. We obtain - using the formulas (4.1-3) - that in
$\NIL$-space for the rotation group $\Gamma_p(q)=\mathbf{pq2_1}$ the metric data of the godesic ball
arrangements $\cB_p(q)$ are the following:
\begin{theorem}
If the system of equation (4.4) holds then the maximal radii and densities of the optimal ball packings are the following:
\begin{itemize}
		\item If $(p,q)=(3,6)$, then $\delta_p(q) \approx 0.2593$, with $R_{opt}(p,q) \approx 0.7389$,
		\item If $(p,q)=(4,4)$, then $\delta_p(q) \approx 0.6512$, with $R_{opt}(p,q) \approx 1.2154$,
		\item If $(p,q)=(6,3)$, then $\delta_p(q) \approx 0.7272$, with $R_{opt}(p,q) \approx 1.9601$.
\end{itemize}
\end{theorem}
If we vary the parameter $x_p(q)$ in the above cases then the corresponding radius $R_{opt}(p,q)$ and the density $\delta_p(q)$ also change.
The following table shows that probably the $\cB_p(q)$ ball packings with maximal kissing numbers provide the optimal ball packing densities.
\begin{center}
  \begin{tabular}{ | c | c | c | c | c |}
    \hline
    $(p,q)$ & Radius & Prism volume & Density & Kissing number \\ \hline \hline
	(3,6) & 0.5876 & 4.1446 & 0.2063 & 2 \\ \cline{2-4}
	 & 0.6392 & 4.9032 & 0.2246 & 2 \\ \cline{2-4}
	 & 0.6929 & 5.7616 & 0.2438 & 2 \\ \cline{2-4}
	 & {\bf 0.7389} & {\bf 6.5517} & {\bf 0.2593} & 8 \\ \cline{2-4}
	 & 0.7787 & 7.8111 & 0.2558 & 6 \\ \cline{2-4}
	 & 0.8132 & 9.0201 & 0.2525 & 6 \\ \cline{2-4}
	 & 0.8481 & 10.3641 & 0.2495 & 6 \\ \hline \hline
	(4,4) & 0.9927 & 7.8849 & 0.5283 & 2\\ \cline{2-4}
	 & 1.0644 & 9.0650 & 0.5678 & 2 \\ \cline{2-4}
	 & 1.1386 & 10.3729 & 0.6090 & 2 \\ \cline{2-4}
	 & {\bf 1.2154} & {\bf 11.8175} & {\bf 0.6512} & 10 \\ \cline{2-4}
	 & 1.2594 & 13.4079 & 0.6404 & 8 \\ \cline{2-4}
	 & 1.3036 & 15.1538 & 0.6295 & 8 \\ \cline{2-4}
	 & 1.3480 & 17.0647 & 0.6194 & 8 \\ \hline \hline
	(6,3) & 1.6934 & 34.4141 & 0.6190 & 2 \\ \cline{2-4}
	 & 1.7801 & 38.0287 & 0.6537 & 2 \\ \cline{2-4}
	 & 1.8690 & 41.9209 & 0.6897 & 2 \\ \cline{2-4}
	 & {\bf 1.9601} & {\bf 46.1044} & {\bf 0.7272} & 14 \\ \cline{2-4}
	 & 2.0087 & 50.5935 & 0.7153 & 12 \\ \cline{2-4}
	 & 2.0573 & 55.4028 & 0.7038 & 12 \\ \cline{2-4}
	 & 2.1059 & 60.5470 & 0.6929 & 12 \\
    \hline
  \end{tabular}
\end{center}

Therefore, we can formulate by the above results the following conjecture:
\begin{conjecture}\label{conj:ball}
The ball arrangements $\cB_p(q)$ provide the densest ball packing arrangements related to $\Gamma_p(q)=\mathbf{pq2_1}$ $\NIL$ isometry group
with parameters $(p,q)=(4,4),(3,6),(6,3)$.
\end{conjecture}
\begin{remark}
The optimal ball packing in the case of $(p,q)=(6,3)$
has a kissing number of $14$, which is greater than the maximal kissing number $12$
in the Euclidean $3$-dimensional space.
In fact, this is the second ball packing arrangement in $\NIL$ that has this high of a kissing number (see \cite{Sz07}).
\end{remark}
\begin{figure}[htbp]\label{pic:prism4}
\centering
\includegraphics[width=180pt]{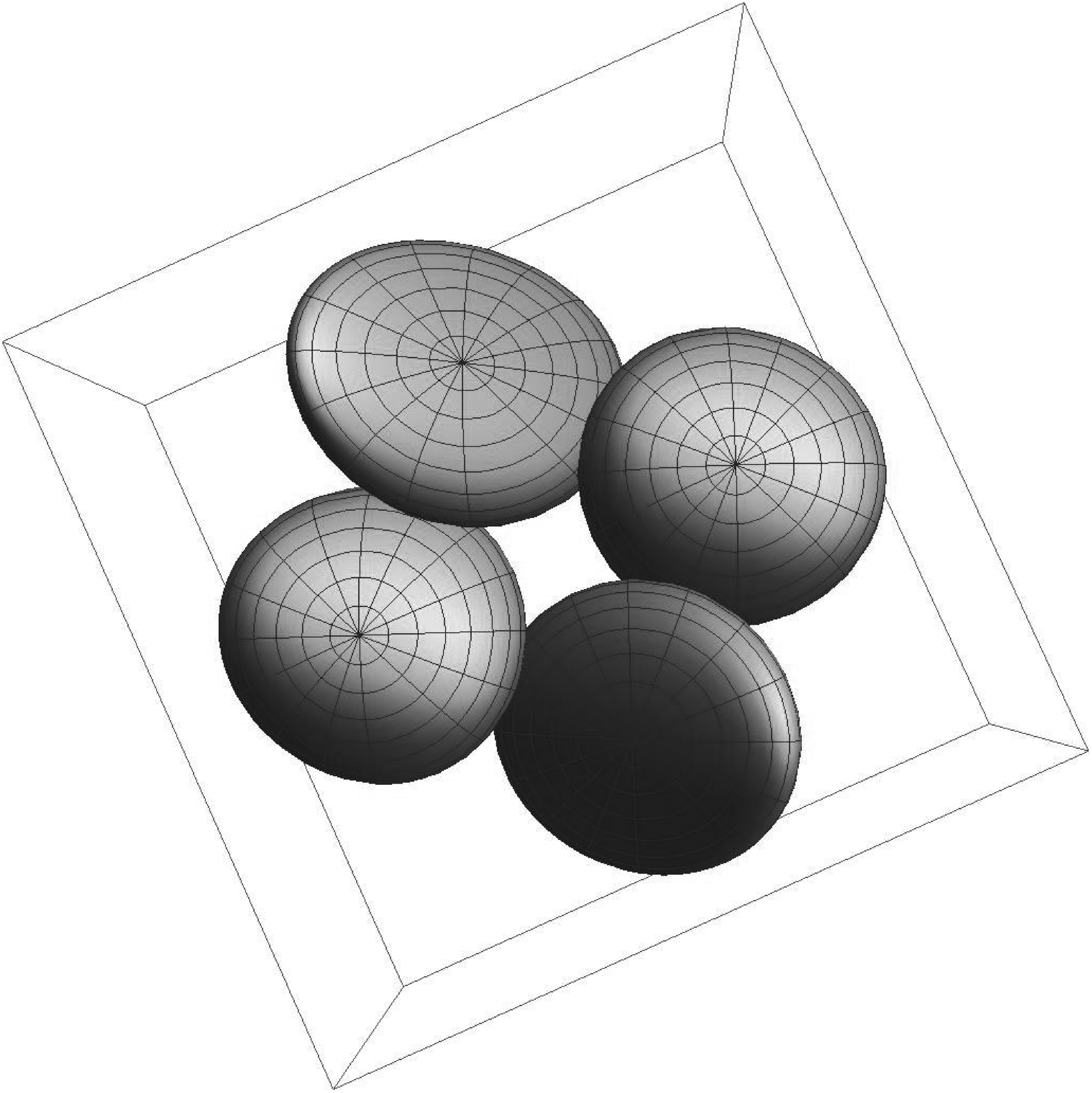}\includegraphics[width=180pt]{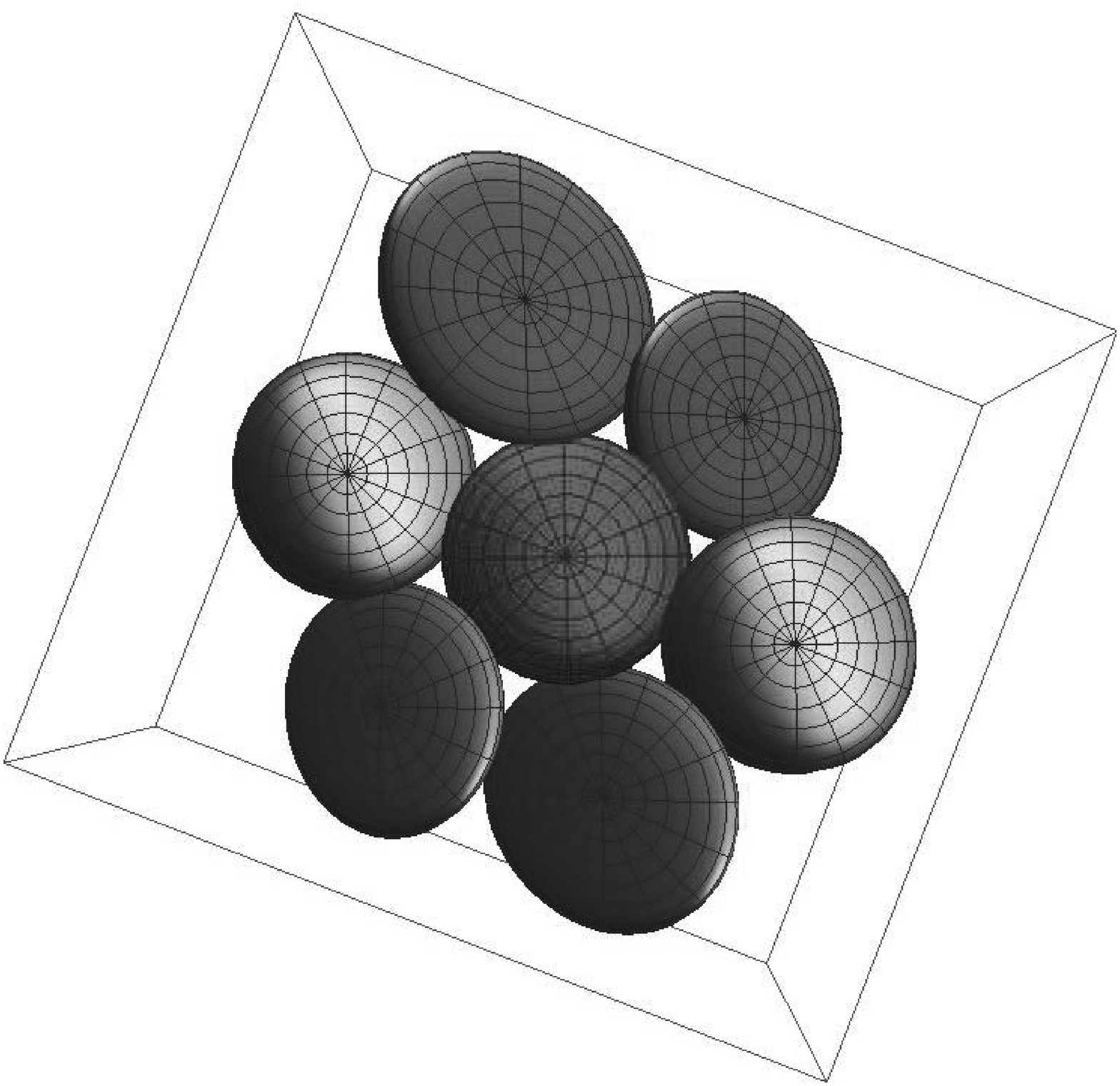}
\caption{Some balls of the optimal ball arrangements for the square and hexagonal
tilings shown from the direction of the $z$-axis.}
\label{pic:optballs2}
\end{figure}
%
%%%%%%%%%%%%%%%%%%%%%%%%%%%%%%%%%%%%%%%%%%%%%%%%%%%%%%%%%%%%%%%%%%%%%%%%%%%%%%%%%%%%
%

\end{document}